\documentclass[12pt]{amsart} 
\usepackage[nice]{nicefrac}
\usepackage{amsmath,amssymb,amsbsy,amsfonts,amsthm,latexsym, amsopn,amstext,amsxtra,euscript,amscd}
\usepackage{delarray}
\usepackage{mathrsfs}
\usepackage{setspace}

\date{}

\renewcommand\overline\bar

\theoremstyle{definition}

\theoremstyle{remark}

\def\ve{\varepsilon}

\def\mod{\text{mod\,}}
\newtheorem*{theorem}{\bf Theorem}
\begin{document}
\title
{Prescribing the binary digits of primes}
\author
{Jean Bourgain}
\begin{address}
{School of Mathematics\\
Institute for Advanced Study\\
 Princeton, NJ 08540}
\end{address}
\thanks
{This research was partially supported by NSF grants DMS-0808042 and DMS-0835373}
\maketitle
\abstract
We present a new result on counting primes $p<N=2^n$ for which $r$ (arbitrarily placed) digits in the binary expansion of $p$ are specified.
Compared with earlier work of Harman and Katai, the restriction on $r$ is relaxed to $r< c\Big(\frac n{\log n}\Big)^{4/7}$. This condition results from the
estimates of Gallagher and Iwaniec on zero-free regions of $L$-functions with `powerful' conductor.
\endabstract

\section*
{(0).  \ Summary}

This work is motivated by the paper [H-K] on the problem described in the title.  We prove the following
\bigskip

\begin{theorem}
{\it Let $N=2^n$, $n$ large enough, and $A\subset\{1, \ldots, n-1\}$ such that
$$
r=|A|< c\left(\frac n{\log n}\right)^{4/7}.\eqno{(0.1)}
$$
Then, considering binary expansions $x= \sum\nolimits_{j< n} x_j 2^j$ $(x_0=1$ and $x_j=0, 1$ for $1\leq j<n)$
and assignments $\alpha_j$ for $j\in A$, we have
$$
|\{ p<N; \text { for } j\in A, \text { the $j$-digit of $p$ equals $x_j$}\}|\sim  2^{-r} \frac N{\log N}.\eqno{(0.2)}
$$}

A few comments.  Statement (0.2) can also be formulated as an asymptotic
formula.
Next, our result has an analogue for $q$-ary expansions (with essentially the same proof) and we restricted ourselves to $q=2$ only for
simplicity.
\end{theorem}

The paper [H-K] establishes a similar result under the assumption
$$
|A|< c \frac {\sqrt n}{\log n}.\eqno{(0.3)}
$$

A key ingredient in [H-K] and also here is the use of zero-free regions for $L$-functions $L(s, \chi)$.
In particular, the restriction (0.1) follows from
 results of Gallagher and Iwaniec ([G], [I]) that provide the zero-free region $1-\sigma < c\frac 1{(\log q
T. \log\log q T)^{3/4}}$, $|\gamma|< T$ where $\rho =\sigma+i\gamma$, for special moduli $q$ that are powers of a fixed integer (here $q=2^j$).
This region is larger than what's available in the general case.
Note that we dismiss here a possible Siegel zero, which is not a concern for $q$ specified as above.
Although the Gallagher-Iwaniec theorem is one of the elements in the [H-K] argument, its full potential was unfortunately not exploited.
As in [H-K], we use the circle method but with a smaller minor arcs region, leading to more saving for that contribution.
As one expects, the major arcs analysis is more involved.

Note that there are two multiplicative character sums entering the discussion (after conversion of the exponential sums in the circle method)
$$
\sum_{n<N} \Lambda(n)  \chi(n)\eqno {(0.4)}
$$
and
$$
\sum_{n<N} f(n)  \chi(n)\eqno{(0.5)}
$$
where $f=1_{[x_j=\alpha_j \text { for } j\in A]}$ or a suitably modified version.

The main additional idea in this paper may be roughly explained as follows.
In general, again dismissing Siegel zeros, $L(s,  \chi)$ has a zero-free region $1-\sigma < \frac c{\log qT}, \hfill\break
|\gamma|~<~T$.
But, as implied by the usual density estimates, for `most' characters, this zero-free region is much larger leading to better bounds on (0.4).
For the remaining `bad' characters, which are few (including possible Siegel zeros), we seek non-trivial estimates on the sum (0.5).
This is clearly reasonable, in view of the additive structure of $f$.
Perhaps not surprisingly, it turns out that the only situation that escapes this analysis (that in fact can be carried out as long $r< n^{\frac 23-\ve})$ are
precisely moduli of the form $2^j$.
Thus at the end, it is the size of the [G], [I] zero-free region that dictates the restriction on $|A|$.

\section {Minor Arcs Contribution}

Let $N=2^n$.
Write
$$
\sum_{k\leq N} \Lambda(k) f(k)=
\int^1_0 S(\alpha) \overline S_f (\alpha) d\alpha\eqno{(1.1)}
$$
denoting
$$
S(\alpha)= \sum\Lambda(k) e(k\alpha)\eqno{(1.2)}
$$
and
$$
S_f(\alpha)= \sum f(k) e(k\alpha).\eqno{(1.3)}
$$
We assume $f(k)=0$ for $k$ even, since obviously $k\equiv 1$ (\mod 2) is a necessary condition.

We fix a parameter $B=B(n)$ which will be specified in \S6.
At this point, let us just say that $\log B\sim n^{4/7}$, up to logarithmic factors.

The major arcs are defined by
$$
\mathcal M(q, a)= \left[\left|\alpha -\frac aq\right|<\frac B{qN}\right] \ \text { where } \  q< B.\eqno{(1.4)}
$$
Given $\alpha$, there is $q<\frac NB$ such that
$$
\left|\alpha-\frac aq\right| <\frac B{qN}<\frac 1{q^2}.
$$

From Vinogradov's estimate  (Theorem 13.6 in [I-K])
$$
|S(\alpha)| < \left(q^{\frac 12} N^{\frac 12}+q^{-\frac 12}N+N^{\frac 45}\right)(\log N)^3\\
$$
$$
\qquad \, \ll \left(\frac N{\sqrt B}+\frac N{\sqrt q}+N^{4/5}\right) (\log N)^3.
\eqno{(1.5)}
$$
Hence if $q\geq B$,
$$
|S(\alpha)|\ll \frac N{\sqrt B} (\log N)^3.\eqno{(1.6)}
$$
Thus the minor arcs contribution in (1.1) is at most
$$
\ll\frac N{\sqrt B}(\log N)^3 \Vert S_f\Vert_1.\eqno{(1.7)}
$$
Write $k=\sum_{0\leq j<n} k_j 2^j$ with $k_j=0, 1$. 
Given $A\subset \{1, \ldots, n-1 \}$ and $\alpha_j=0$ or $1$ for $j\in A$, we define
$$
f(x) =\prod_{j\in A} h_{\alpha_j} \left(\frac x{2^{j+1}}\right)\eqno{(1.8)}
$$
where $0\leq h_0 \leq 1$ is a smooth 1-periodic function satisfying 
$$
\begin{cases} h_0(t)=1 &\text { if } \frac 1n\leq t\leq \frac 12-\frac 1n\\ h_0(t)=0 &\text { if } \frac 12\leq t\leq 1\end{cases}
$$
($h_1$ is defined similarly).
Thus
$$
\Vert \hat h_\alpha\Vert_1 < C \log n\eqno{(1.9)}
$$
and, by approximation, we can assume that say
$$
\text{supp\,} \hat h_\alpha \subset [-n^3, n^3]\eqno {(1.9')}
$$
the additional error term (for suitably chosen $h_\alpha$) being at most $e^{-n} <\frac 1N$ and hence may be ignored.

Note that in the formulation of the main theorem, we did not specify the error term as we do not intend to put emphasis on this aspect.
Thus the function $f$ introduced in (1.8) suffices for our purpose, while establishing an asymptotic formula with specified error term would require an
approximation of the functions $1_{[0, \frac 12[}$ and $1_{[\frac 12, 1[}$ which is better than given by $h_0$ and $h_1$.

In particular $f$ satisfies
$$
\Vert S_f\Vert_1 <(C\log n)^{|A|}.\eqno{(1.10)}
$$
Substitution in (1.7) gives then the bound
$$
\frac N{\sqrt B}(C\log n)^{|A|} n^3.\eqno{(1.11)}
$$

\section{Major Arcs Analysis (I)}

Next, we analyze the major arcs contributions $(q< B)$
$$
\sum_{(a, q)=1} \ \int\limits_{|\alpha-\frac aq|<\frac B{qN}} S(\alpha) \overline S_f(\alpha) d\alpha.\eqno{(2.1)}
$$
Write $\alpha =\frac aq +\beta$. Defining

$$
\tau(\chi) =\sum^q_{m=1}  \chi(m) e_q(m)
$$
we have (see [D], p. 147)
$$
S(\alpha) =\frac 1{\phi(q)} \sum_{\chi} \tau (\overline{\chi}) \chi(a) \left[\sum_{k\leq N} \chi(k) \Lambda (k)
e(k\beta)\right]+O\big((\log N)^2\big)\eqno{(2.2)}
$$
where by partial summation
$$
\left|\sum_{k\leq N} \chi(k) \Lambda (k) e(k\beta)\right|\leq (1+N|\beta|) \max_{u<N}|\psi (u, \chi)|
$$
$$
\qquad\qquad <\frac Bq \, \max_{u<N}|\psi (u,  \chi)|.\eqno{(2.3)}
$$
We subdivide the multiplicative characters $\chi$ in classes $\mathcal G$ and $\mathcal B$ (to be specified), depending on the zero set
of $L(s,  \chi)$.

Thus
$$
S(\alpha)= \frac {\mu(q)}{\phi (q)} \left[\sum_{k\leq N} \Lambda (k) e(k\beta)\right] +(2.4) + (2.5)
$$
where the first term comes from the contribution of the principal character $\chi_0$ and
$$
(2.4) =\frac 1{\phi(q)} \sum_{\substack {\chi(\mod q)\\ \chi\in\mathcal G}} \tau(\overline {\chi})\chi(a) \left[\sum_{k\leq N}
\chi(k) \Lambda (k) e(k\beta)\right]
$$
and
$$
(2.5) =\frac 1{\phi(q)} \sum_{\substack{\chi(\mod q) \\ \chi\in \mathcal B}} \cdots \qquad\qquad\qquad\qquad\qquad\qquad
$$

Hence
$$
|(2.4)|\lesssim \frac B{\sqrt q} \left[\max_{\substack{u<N \\ \chi(\mod q), \chi\in\mathcal G}}  |\psi(u,  \chi)|\right]\eqno{(2.6)}
$$
by (2.3).
The contribution of (2.6) in (2.1) is at most
$$
(2.6). \quad  \int\limits_{\bigcup_{(a, q)=1}\mathcal M(q, a)} |S_f(\alpha) |  d\alpha.
$$
Summing over $q<B$ and using the bound (1.10) gives the estimate
$$
B(C\log n)^{|A|} \Bigg\{ \max_{\substack{u<N \\  \chi\in\mathcal G}} |\psi (u,  \chi)|\Bigg\}\eqno{(2.7)}
$$
where in (2.7) the conductor of $\chi$ is at most $B$.

\section{Major Arcs Analysis (II)}

Next we treat $\chi=\chi_0$ and $ \chi\in \mathcal B$.

Assume $\chi$ is induced by $\chi_1$ which is primitive $(\mod q_1)$, $q_1|q$.
Then from [D], p. 67
$$
\tau (\overline{\chi})= \mu\left(\frac q{q_1}\right) \overline{\chi}_1 \left(\frac q{q_1}\right) \tau \left(\overline{\chi}_1\right)\eqno{(3.1)}
$$
which vanishes, unless $q_2 =\frac q{q_1}$ is square free with $(q_1, q_2)=1$.

The contribution of $\chi$ in (2.1) equals
$$
\frac{\tau(\overline{\chi})}{\phi (q)} \int\limits_{|\beta|<\frac B{qN}}\left[\sum_{k\leq N} \chi(k) \Lambda(k) e(k\beta)\right]
\left[\sum_{k<N} f(k) \left(\sum^q_{a=1} \chi(a) e_q(-ak)\right) e(-k\beta)\right] d\beta.\eqno{(3.2)}
$$
We have
$$
\begin{aligned}
\sum^q_{a=1} e_q(ak) \chi(a) &=\sum _{(a, q)=1} e_q(ak) \chi_1(a)
\\[6 pt]
&=\Bigg[\sum_{(a_1, q_1)=1} e_{q_1} (a_1k) \chi_1 (a_1)\Bigg] \Bigg[\sum_{(a_2, q_2)=1} e_{q_2} (a_2k)\Bigg]
\end{aligned}
$$
$$
= \chi_1(k) \tau(\chi_1) c_{q_2} (k)\qquad \qquad \eqno{(3.3)}
$$
(cf. [D], p. 65) and
$$
c_{q_2}(k) =\frac {\mu \left(\frac{q_2}{(q_2, k)}\right) \phi(q_2)}{\phi\left(\frac{q_2}{(q_2, k)}\right)}.\eqno{(3.4)}
$$
(cf. [D], p. 149).

From (3.1), (3.3), (3.4)
$$
\frac {\tau(\overline{\chi})}{\phi(q)}\left[ \sum^q_{a=1}  \chi(a) e_q(-ak)\right] =\frac {|\tau({\chi_1})|^2}{\phi(q_1)} \ 
\frac {\bar{\chi}_1 (q_2)}{\phi\left(\frac {q_2}{(q_2, k)}\right)}\,  \mu\big(( q_2, k)\big) \mathcal
\chi_1(k)\qquad\qquad\qquad
$$
$$
\qquad\qquad =\frac {q_1}{\phi(q_1)} \ \frac {\overline{\chi}_1 (q_2)}{\phi(\frac{q_2}{(q_2, k)})} \, \mu\big((q_2, k)\big)  \chi_1(k)
\eqno{(3.5)}
$$
(cf. [D], p. 66).

Returning to (3.2), rather than integrating
in $\beta$ over the interval $|\beta|<\frac B{qN}$, we introduce a weight function
$$
w\Big(\frac{qN}{B} \beta\Big)
$$
where $0\leq w\leq 1$ is a smooth bumpfunction on $\mathbb R$ such that $w=1$ on $[-1, 1]$, supp $w\subset [-2, 2]$ and
$$
|\hat w(y)|< C e^{-|y|^{1/2}}.
$$
\big(Note that this operation creates in (2.1) an error term that is captured by the minor arcs contribution (1.7)\big).

Hence, substituting (3.5), (3.2) becomes
$$
\frac {q_1}{\phi(q_1)} \overline {\mathcal X}_1 (q_2)\frac B{qN} \sum_{k_1, k_2<N} \hat w\Big(\frac B{qN} (k_1-k_2)\Big)\mathcal X(k_1)\Lambda (k_1) f(k_2)\frac
{u\big((q_2, k_2)\big)}{\phi\big(\frac {q_2}{(q_2, k_2)}\big)} \mathcal X_1(k_2)\eqno{(3.6)}
$$
and we observe that by our assumption on $\hat w$ the $k_1, k_2$ summation in (3.6) is restricted to $|k_1-k_2|< \frac {qN}B n^3$, up to a negligible error.

We first examine the contribution of the principal characters.

For $\mathcal X=\mathcal X_0, q_1=1$ and (3.6) becomes
$$
\frac B{qN} \sum_{k_1, k_2<N} \hat w \Big(\frac B{qN} (k_1-k_2)\Big) \Lambda(k_1) f(k_2) \frac {\mu\big(( q, k_2)\big)}{\phi\big(\frac q{(q, k_2)}\big)}.\eqno{(3.7)}
$$
Recall the distributional property of the primes ([I-K], Corollary 8.30)
$$
\sum_{k\leq x}\Lambda(k) =x+O\big(x\exp \big(-c(\log x)^{3/5}(\log\log x)^{-\frac 15}\big)\big)\eqno{(3.8)}
$$
and also our assumption $\log B\ll n^{4/7}$.
Performing the $k_1$-summation in (3.7), partial summation together with (3.8) give
$$
\begin{aligned}
&\frac B{qN} \sum_{k_1, k_2} \hat w \Big(\frac B{qN} (k_1-k_2)\Big) f(k_2) \frac {\mu\big((q, k_2)\big)}{\phi\big(\frac q{(q, k_2)}\big)}
+O\big(N\exp \big(-c (\log N)^{3/5}(\log\log N)^{-\frac 15}\big)\big)\\
& =\sum_k f(k)\frac {\mu\big((q, k)\big)}{\phi\big(\frac q{(q, k)}\big)}+ O\big(N\exp \big(-c(\log N)^{3/5}(\log\log N)^{-\frac 15}\big)\big)
\end{aligned}
\eqno{(3.9)}
$$

Let $\kappa(q)$ be a function satisfying the following
\null\vskip.2 true in

\noindent
{\bf Assumption A.} {\sl Let $q_0<B$ be odd and square free. Then
$$
\sum_{q_0|k, k<N} f(k) =\mathbb E[f] \frac N{q_0} +O\Big(\kappa(q_0) \frac N{q_0} 2^{-r}\Big)\eqno{(3.10)}
$$
where $\mathbb E[f]$ denotes the normalized average.}
\null\vskip.3 true in

\noindent
{\bf Remark.} Since the function $\kappa(q)$ at this stage remains unspecified and its properties will be established later (in \S4), we 
thought this presentation more desirable than making a more technical statement here (cf. (4.26)).

Assuming $q$ square-free (sf)  and  odd, the first term of (3.9) equals
$$
\begin{aligned}
&\sum_{q'|q} \, \frac {\mu(q')}{\phi\big(\frac q{q'})} \left[\sum_{(q, k)=q'} f(k)\right]=\\[6pt]
&\sum_{q'|q} \, \frac{\phi(q')}{\phi(q)} \mu(q') \sum_{q''|\frac q{q'}} \mu(q'')\left[\sum_{\substack {q'q''|k\\ k<N}} f(k)\right]
\end{aligned}
$$
and substituting (3.10) we obtain
$$
N \, \mathbb E [f] \, \sum_{q'|q} \, \frac{\mu(q')}{\phi(q/q')} \, \sum_{q''|\frac q{q'}} \, \frac {\mu(q'')}{q'q''}\eqno{(3.11)}
$$
$$
+O\left( N\, 2^{-r} \sum_{q'|q} \, \frac {\phi(q')}{\phi(q)} \, \sum_{q''|\frac q{q'}} \, \frac {\kappa (q'q'')}{q'q''}\right).\eqno{(3.12)}
$$

Next
$$
\begin{aligned}
(3.11) &= N\, \mathbb E[f] \, \sum_{q'|q} \, \frac {\mu(q')}{\phi (q/q')} \, \frac 1{q'} \prod_{p|\frac q{q'}} \, \left(1-\frac 1p\right)
\\[6pt]
&=N \,\mathbb E[f] \, \sum_{q'|q} \, \frac {\mu(q')}{\phi (q/q')} \, \frac 1q \, \phi (q/q')
\end{aligned}
$$
$$
=\begin{cases} N\, \mathbb E[f] &\text { if }  \ q=1\\ 0&\text { otherwise.}\end{cases}\qquad \eqno{(3.13)}
$$

Summing (3.12) over $q<B$ sf and odd, we have the estimate (setting $q_1=q'q''$)
$$
2^{-r} N \sum_{\substack{q_1< B\\ q_1 \, sf}} \, \frac {\kappa(q_1)}{q_1} \left[\sum_{\substack{ q''|q_1, (q_1, q_2)=1\\ q_2< B \, sf}}
\, \frac 1{\phi(q'')\phi (q_2)}\right]
$$
$$
<2^{-r} N(\log B)^2 \left[\sum_{\substack{q<B\\ q\, sf, \text { odd}}}\ \frac {\kappa(q)}{q}\right].\eqno{(3.14)}
$$
For $q$ $sf$ and even, set $q=2q_1$ and note that the first term of (3.9) equals
$$
\sum_{q'|q_1} \ \frac{\mu(q')}{\phi\big(q_1/q')} \ \left[\sum_{(q_1, k)=q'} f(k)\right]
$$
and we proceed similarly as above with the same conclusion and $q$ replaced by $q_1$.
From the preceding, the contribution of the principal characters equals
$$
2\mathbb E[f] N + CN2^{-r}(\log N)^2 \left[\sum_{\substack {q<B\\ q \, sf}} \ \frac {\kappa(q)}{q}\right]+ BN\exp \big(-c(\log N)^{3/5} (\log\log
N)^{-\frac 15}\big).
\eqno{(3.15)}
$$
Consider next $\chi_1\in \mathcal B, \chi_1$ primitive.

Returning to (3.6), estimate by
$$
\frac {q_1}{\phi(q_1)} \ \frac B{qN} \sum_{k_1<N} \Big|\sum_k \hat w\Big(\frac B{qN} (k_1-k)\Big) f(k)\frac {\mu\big((q_2, k)\big)}{\phi\big(\frac {q_2}
{(q_2, k)}\big)} \mathcal X_2(k)\Big|\qquad\qquad
$$
$$
\leq \frac {q_1}{\phi(q_1)}\ \frac Bq\sum_{q_2' | q_2}\frac {\phi(q_2')}{\phi(q_2)} \max_{k_1} \Big|\sum_{k, (k, q_2)=q_2'} \hat w\Big(\frac B{q_N}
(k_1-k)\Big) f(k)\mathcal X_1(k)\Big|
\eqno{(3.16)}
$$
\medskip

We introduce another function $\alpha(q_1, q_0)$ satisfying
\bigskip

\noindent
{\bf Assumption B.}

{\sl	Given $q_0, q_1 <B, (q_0, q_1)= 1$ with $q_0$ $sf$ and odd, $ \chi_1(\mod q_1)$ primitive,
$$
\left|\sum_{k\in J, q_0|k} \ f(k)  \chi_1(k)\right| < 2^{-r} \alpha(q_1, q_0)\, \frac{|J|}{q_0}\eqno{(3.17)}
$$
holds, whenever $J\subset[1, N]$ is an interval of size $\sim\frac NB$.}

The function $\alpha(q_1, q_0)$ which enters in (3.17) will be studied in \S4, \S5 and in particular can be chosen to satisfy (5.11).
Partial summation using the decay estimate on $\hat w$ and (3.17) implies the following bound on (3.16)
$$
2^{-r}N (\log N)^3\ \frac {q_1}{\phi(q_1)} \ \sum_{q_2'|q_2} \ \frac {\phi(q_2')}{\phi(q_2)} \ \sum_{q_2''|\frac {q_2}{q_2'}}
\ \frac {\alpha(q_1, q_2' q_2'')}{q_2' q_2''}.\eqno{(3.18)}
$$
Summation of (3.18) over $sf$ $q_2<B, (q_1, q_2)=1$ gives the following bound for the $\mathcal X_1$-contribution
$$
2^{-r} N (\log N)^3 \frac {q_1}{\phi(q_1)} \ \sum_{\substack{q_3<B\\ sf, \text { odd}}} \ \frac{\alpha (q_1, q_3)}{q_3} \ \sum_{\substack
{q_2''|q_3 \\ q_2'''<B, (q_2''', q_2'')=1}} \ \frac 1{\phi(q_2'')\phi(q_2''')}\\[6pt]
$$
$$
<2^{-r} N (\log N)^5 \ \frac {q_1}{\phi(q_1)} \, \left\{\sum_{\substack{q_3< B\\ sf \text { odd}}} \ \frac {\alpha(q_1, q_3)}{q_3}
\right\}.\eqno{(3.19)}
$$
The necessary information on the functions $\kappa(q)$ and $\alpha(q_1, q_0)$ needed to bound (3.15) and (3.19) will be established in
the next sections \S4, \S5.

\section{Further Estimates (I)}

Consider condition (3.10) and the expression (3.14).

Let
$$
A=\{j_1< j_2< \cdots< j_r\} \subset \{1, \ldots, n\}.
$$
Let $1< J_1< J_2<n$ (to be specified) and write
$$
x=\sum_{j<J_1} x_j 2^j+\sum_{J_1\leq j<J_2} x_j 2^j+\sum_{J_2\leq j<n} x_j2^j =u+v+w.
$$
We evaluate (assuming $q_0$ odd)
$$
\sum_{q_0|x} f(x) =\sum_{u, w} \left[\sum_{q_0|u+v+w} f(u+v+w)\right]
$$
by evaluating the inner sum in $v$ with $u, w$ fixed.

Clearly
$$
f(x)=\prod_{j\in A} h\left(\frac x{2^{j+1}}\right)=
\qquad\qquad\qquad\qquad
$$
$$
\prod_{\substack {j\in A\\ j\leq J_1}} h\left(\frac u{2^{j+1}}\right). \prod_{\substack{j\in A\\ J_1<j\leq J_2}} \left[ h\left(\frac v{2^{j+1}}\right)+O(\Vert
h\Vert_{C^1} 2^{J_1-j})\right].
\prod_{\substack{j\in A\\ J_2<j\leq n}}\left[h \left(\frac w{2^{j+1}}\right)+ O(\Vert h\Vert_{C^1} 2^{J_2-j})\right].\eqno{(4.1)}
$$
Fix $\Delta J = J_2-J_1>10\log B$ and take $J_1, J_2$ satisfying
$$
\text{dist } (\{J_1, J_2\}, A)>\frac{n}{10r}\eqno{(4.2)}
$$
and
$$
|A\cap [J_1, J_2]|<\frac {10 r\Delta J}{n}.\eqno{(4.3)}
$$
Since $\Vert h\Vert_{C^1}\lesssim n^2$, assuming
$$
r<10^{-6} \frac n{\log n}\eqno{(4.4)}
$$
we have
$$
f_-\leq f\leq f_+
$$
with
$$
f_\pm (x) =\prod_{\substack{j\in A\\ j\leq J_1}} h\left(\frac u{2^{j+1}}\right) .\prod_{\substack{j\in A\\ J_1<j\leq J_2}}
h_\pm \left(\frac v{2^{j+1}}\right). \prod_{\substack{j\in A\\ j>J_2}} h_\pm \left(\frac w{2^{j+1}}\right)
$$
and
$$
h_\pm = h\pm 2^{-\frac n{11r}}.
$$
Note that
$$
\mathbb E[f_+]- \mathbb E[f_-]\lesssim r2^{-\frac n{11 r}} 2^{-r} < 2^{-\frac n{12r}}2^{-r}.
$$

Obviously
$$
\sum_{q_0|x} f_-(x)\leq \sum_{q_0|x} f(x) \leq \sum_{q_0|x} f_+(x)
$$
and hence
$$\left|\sum_{q_0|x} f(x)-\frac N{q_0} \mathbb E[f]\right|<
\max_{\pm} \left|\sum_{q_0|x} f_\pm (x) -\frac N{q_0} \mathbb E[f_\pm]\right|+2^{-\frac n{12r}} \frac N{q_0} 2^{-r}.
$$
Assume that for fixed $u, w$
$$
\left|\sum_{\substack {v\equiv -u-w\\ (\mod q_0)}} \,  \prod_{\substack {j\in A\\ J_1< j\leq J_2}} h_+ \left(\frac v{2^{j+1}}\right)-\frac {2^{\Delta
J}}{q_0} \mathbb E\left[ \prod_{\substack{j\in A\\ J_1<j\leq J_2}} h_+ \left(\frac .{2^{j+1-J_1}}\right)\right]\right| <
$$
$$
\kappa (q_0) \frac {2^{\Delta J}}{q_0} 2^{-|A\cap [J_1, J_2]|}.\eqno{(4.5)}
$$
It will follow that
$$
\left|\sum_{q_0|x} f_+(x) -\frac N{q_0} \mathbb E[f_+]\right|< \kappa (q_0) 2^{-r} \frac N{q_0}
$$
(and similarly for $f_-$).

Denote
$$
 A'=A\cap ]J_1, J_2]= J_1-1  +\{j_1' <\cdots< j_{r_1}'\} \ \text { where } \ r_1< \frac {10r\Delta J}n.
$$

Expanding as a Fourier series, we get
$$
\prod_{j\in A'} h_+\left(\frac v{2^{j+1}}\right) -\mathbb E[\cdots] = \sideset{}{^*}\sum_{\{b_j\}_{j\in A'}}\, \left[\prod_{j\in A'} \widehat h_+(b_j)\right]
e\left(\sum_{j\in A'} \, \frac {b_j}{2^{j+1}} v\right)\eqno{(4.6)}
$$
where $\sum^*$ refers to those $\{b_j\}_{j\in A'}$ with $\sum_{j\in A'}$ $\frac {b_j}{2^{j+1 -J_1}}\not= 0 (\mod 1)$.

Recall that $|b_j|<n^3$ and
$$
\sum_{\{b_j\}} \ \prod_{j\in A'} |\widehat h_+(b_j)|<(C\log n)^{r_1}.\eqno{(4.7)}
$$
Clearly, the left side of (4.5) is at most
$$
\frac {2^{\Delta J}}{q_0} \sideset{}{^*}\sum_{\{b_j\}} \ \prod_{j\in A'} |\widehat h_+(b_j)|\left\{1_{\big
[\Vert q_0 (\sum_{j\in A'}\frac {b_j}{2^{j+1 -J_1}})\Vert<\frac {n^2q_0}{2^{\Delta J}}\big]} +2^{-n}\right\}
$$
and we may take in (4.5)
$$
\kappa(q_0) =2^{r_1}\sideset{}{^*}\sum_{\{b_j\}} \ \prod_{j\in A'} |\widehat h_+ (b_j)|. 
1_{\big[\Vert q_0 \sum_{j\in A'} \frac {b_j}{2^{j+1-J_1}}\Vert<\frac {n^2q_0}{2^{\Delta
J}}\big]}.\eqno{(4.8)}
$$

Denote $\beta_s =b_{j_s'} (1\leq s\leq r_1)$.
Fix $\{\beta_s\}$ such that
$$
\sum^{r_1}_{s=1} \, \frac {\beta_s}{2^{j_s'}} \not= 0(\mod 1)\eqno{(4.9)}
$$
and estimate for given $Q<B$
$$
\sum_{\substack {q\sim Q\\ q \text { odd}}} 1_{\big[\Vert q(\sum_s \frac {\beta_s}{2^{j_s'}})\Vert < \frac {n^2Q}{2^{\Delta J}}\big]}.\eqno{(4.10)}
$$
Since $q$ is odd, $q\left(\sum_s\frac {\beta_s}{2^{j_s'}}\right) =\frac b{2^{j_{r_1}'}} \not= 0 (\mod 1)$ and, if $(4.10)\not=0$,\ necessarily
\hfill\break $2^{-j_{r_1}'}< n^2 Q 2^{-\Delta J}$, $2^{-\frac n{10r}} \overset {(4.2)}< 2^{J_2 -J_1-j_{r_1}'}< n^2 Q$, implying
$$
Q> 2^{\frac n{11r}}.\eqno{(4.11)}
$$

Let $(4.10)=\frac QK$, $K\leq Q$.
From the pigeon hole principle, there is some $1\leq q_1 \leq K$ and $0\leq a_1 < q_1$ with
$$
\left\Vert \sum \ \frac{\beta_s}{2^{j_s'}} -\frac {a_1}{q_1}\right\Vert < \frac {n^2Q}{2^{\Delta J}}.\eqno{(4.12)}
$$
Note that $q|\sum\frac {\beta_s}{2^{j_s'}}|<4 Q\frac {n^3}{2^{j_1'}}$.
If $2^{j_1'}> 8Bn^3$, it would follow that $\Vert q. \big(\sum \ \frac {\beta_s}{2^{j_s'}}\big)\Vert =q \big|\sum \frac{\beta_s}{2^{j_s'}}\big|<
Q n^2 2^{-\Delta J}$ and $2^{-j_{r_1}'}< 2n^j 2^{-\Delta J}$, contradicting (4.2).

Hence $2^{j_1'}\leq 8 n^3B$.

Let $D\gg 1$ be another parameter with
$$
r_1D<\frac 1{10} \Delta J.\eqno{(4.13)}
$$
Clearly there is some $1\leq s_0 \leq r_1$ such that
$$
j_{s_0}' < j_1' +r_1D\eqno{(4.14)}
$$
and, if $s_0<r_1$,
$$
j_{s_0+1}'-j_{s_0}' > D.\eqno{(4.15)}
$$
Hence
$$
2^{j_{s_0}'} < 8n^3 B2^{r_1 D}.\eqno{(4.16)}
$$
By (4.12), (4.15)
$$
\left\Vert\sum_{s\leq s_0} \frac{\beta_s}{2^{j_s'}} -\frac {a_1}{q_1} \right\Vert < \frac {n^2B}{2^{\Delta J}} +\frac {n^3}{2^{j_{s_0}'+D}}<\frac 1
{2^{j_{s_0}'+\frac 12 D}}\eqno{(4.17)}
$$
assuming
$$
D>10\log n \text { and } \Delta J>10\log B\eqno{(4.18)}
$$
(4.17) implies that either $q_1\geq 2^{\frac 12 D}$, hence
$$
K>2^{\frac 12 D}\eqno{(4.19)}
$$
or
$$
\sum_{s\leq s_0} \ \frac {\beta_s}{2^{j_s'}} =\frac {a_1}{q_1}\qquad (\mod 1).\eqno{(4.20)}
$$
Assuming (4.20), we deduce from (4.10), (4.12) that
$$
\left\Vert qq_1 \left(\sum_{s>s_0} \ \frac {\beta_s}{2^{j_s'}}\right) \right\Vert
< \frac {n^2 q_1 Q}{2^{\Delta J}}
$$
and
$$
\left\Vert\sum_{s>s_0} \ \frac {\beta_s}{2^{j_s'}}\right\Vert <\frac {n^2Q}{2^{\Delta J}}.
$$
By (4.18) and the preceding, $\Vert qq_1 \big(\sum_{s>s_0} \, \frac {\beta_s}{2^{j_s'}}\big)\Vert =qq_1 \Vert\sum_{s>s_0} \, \frac {\beta_s}{2^{j_s'}}\Vert$
so that \hfill\break
$\Vert\sum_{s>s_0} \, \frac {\beta_s}{2^{j_s'}}\Vert<\frac {2n^2}{2^{\Delta J}}$.
Again using that $j_{r_1}' <\Delta J-\frac n{10r}$, it follows that  $\sum_{s>s_0} \, \frac {\beta_s}{2^{j_s'}}=0$ $(\mod 1)$  and
$$
\sum_s \, \frac {\beta_s}{2^{j_s'}}=\frac {a_1}{q_1} \not= 0 (\mod 1)\eqno{(4.21)}
$$
recalling (4.9).

But then from (4.10), since $q_1$ is a power of 2 and $q$ is odd
$$
\frac 1{q_1}\leq\left\Vert q \frac {a_1}{q_1}\right\Vert <\frac {n^2 B}{2^{\Delta J}}
$$
(a contradiction).

Thus (4.19) holds and
$$
(4.10) < 2^{-\frac 12D} Q.\eqno{(4.22)}
$$

Returning to (4.8), (4.7), (4.22) and the bound $r_1<\frac {10r\Delta J}n$ imply that
$$
\sum_{\substack{q\sim Q\\ q\text { odd}}}  \kappa(q) < (C\log n)^{\frac {10r\Delta J}n} \ 2^{-\frac 12 D} Q.\eqno{(4.23)}
$$
Recalling (4.13), take $D= \frac 1{200} \frac nr$.
Taking $\Delta J\sim \log B$ and assuming
$$
\log B< c(\log\log n)^{-1} \, \frac {n^2}{r^2}\eqno{(4.24)}
$$
we obtain
$$
(4.23) < 2^{-10^{-3} \frac nr} Q.\eqno{(4.25)}
$$
Returning to (3.14), it follows in particular from the preceding that
$$
N2^{-r} (\log B)^2 \sum_{\substack {q< B\\ q\, sf, \text { odd}}} \ \frac {\kappa (q)}q<
$$
$$
 N2^{-r}(\log B)^2 (2^{-\frac n{12r}} .\log B+2^{-10^{-3}\frac nr} \log B)<
N2^{-r} 2^{-10^{-4}\frac nr} 
\eqno{(4.26)}
$$
assuming (4.4), (4.24).
\null\vskip.3 true in

Going over our analysis, it is clear that in (4.5) the condition $v\equiv -u-w (\mod q_0)$ can be replaced by $v\equiv a (\mod q_0)$
for arbitrary $a$.
Also, $J_1$ may be taken at least $\frac n2$, so that the element $v\equiv 0 (\mod 2^{[\frac n2]})$.

This variant of the above argument gives the following statement.

Let $q_0<B$ be odd and $m<\frac n2$.  Then
$$
\sum_{\substack{x\equiv a_0(\mod q_0)\\ x\equiv a'(\mod 2^m)}} \ f(x) =\frac 1{q_0} \left[1+0\big(\kappa(q_0)+2^{-\frac n{12r}}\big)\right]
\left[\sum_{\substack{x\equiv a'\\ (\mod 2^m)}} \ f(x)\right].\eqno{(4.27)}
$$
From (4.27), we derive a character sum estimate relevant to (3.17).

Let $q_1 =\tilde q_1 2^m(1< \tilde q_1 \text { odd})$, $q_1 <B$ and $\chi_1=\tilde{\chi}_1 \chi'$ with $\tilde{\chi}_1$ non-principal
$(\mod \tilde q_1)$ and $ \chi'(\mod 2^m)$.
Let $q_3 <B$ be square free, odd and $(q_1, q_3)=1$.

We apply (4.27) with $q_0=\tilde q_1 q_3$. Clearly
$$
\left|\sum_{q_3|x}  \chi_1(x) f(x)\right| \leq \sum_{a'=1}^{2^{m_{-1}}}\Bigg| \sum_{a=1}^{\tilde q_1}\tilde{\chi}_1(a) \Bigg\{
\sum_{\substack{ x\equiv a'(\mod 2^m)\\ x\equiv a(\mod \tilde q_1)\\ q_3|x}}\ f(x)\Bigg\}\Bigg|
$$
$$
\overset{(4.27)}\lesssim \frac 1{q_3} \big(\kappa(q_3 \tilde q_1)+2^{-\frac n{12r}}\big) 2^{-r} N.\eqno{(4.28)}
$$
Similarly, if $J\subset [1, N]$ is an interval of size $\frac NB$
$$
\left|\sum_{x\in J, q_3|x}  \chi_1(x) f(x)\right|\lesssim \frac 1{q_3} \big(\kappa(q_3 \tilde q_1)+ 2^{-\frac n{12r}}\big)
2^{-r}|J|.\eqno{(4.29)}
$$
Hence, in (3.17)
$$
\alpha(q_1, q_3)\lesssim \kappa (q_3\tilde q_1)+2^{-\frac n{12r}}.\eqno{(4.30)}
$$
Recall (4.23), (4.25), it follows that in (3.19)
$$
\begin{aligned}
\sum_{q_3<B} \ \frac{\alpha(q_1, q_3)}{q_3} &\lesssim (\log B) 2^{-\frac n{12r}}+\sum_{q_3<B} \ \frac {\kappa(q_3 \tilde q_1)}{q_3}\\[12pt]
&\lesssim (\log B) 2^{-\frac n{12r}} +\tilde q_1 \sum_{q< B\tilde q_1} \ \frac {\kappa(q)}q
\end{aligned}
$$
$$
\lesssim \tilde q_1 2^{-10^{-3}\frac nr} \log N \eqno{(4.31)}
$$
which is conclusive for $1<\tilde q_1 < 2^{\frac 12 10^{-3}\frac nr}$.

\section{Further Estimates (II)}

Next, we establish a second bound on $\alpha(q_1, q_0)$.

Take $J=[1, N]$; the adaptation for $J$ an interval of size $\frac NB$ is straightforward.

Proceeding as in \S4, take $J_1<J_2$ satisfying (4.2), (4.3) and $\Delta J=J_2=J_1 \sim \log B$.
Write $x=\sum_{j<J_1} x_j2^j+\sum_{J_1\leq j< J_2} x_j 2^j+\sum_{J_2\leq j<n} x_j 2^j =u+v+w$.

Evaluate
$$
\Big|\sum_{q_0|k} f(k) \mathcal X_1(k)\Big|\leq \sum_{u, w} \Big| \sum_{v\equiv -u-w} f(u+v+w)\mathcal X_1 (u+v+w)\Big|\eqno{(5.1)}
$$
by evaluating the inner sum with $u, w$ fixed ($\equiv$ refers to $(\mod q_0)$).

Using the notations from \S4
$$
(5.1)\leq \sum_{u, w} \prod_{j\in A, j<J_1} h\Big(\frac u{2^j}\Big) \prod_{j\in A, j>J_2} h\Big(\frac w{2^j}\Big) \Big|\sum_{v\equiv -u-w}
\prod_{j\in A, J_1< j<J_2} h\Big(\frac v{2^j}\Big) \mathcal X_1 (u+v+w)\Big|\eqno {(5.2)}
$$
$$
+2^{-\frac n{12r}} \sum_{q_0|k} f_+(k).\eqno{(5.3)}
$$
It follows in particular from the analysis in \S4 that
$$
(5.3) < 2^{-\frac n{12r}} 2^{-r}\frac N{q_0}.\eqno{(5.4)}
$$
Denoting $F(v)= \prod_{j\in A, J_1< j<J_2} h\Big(\frac v{2^j}\Big)$, we have \big(cf. (4.7)\big)
$$
\begin{aligned}
\Vert \hat F\Vert_1 &< (C\log n)^{|A\cap [J_1, J_2]|}\\
&<(C\log n)^{10\frac rn\Delta J}\qquad\qquad \text { by } \  (4.3)
\end{aligned}
$$
$$
< B^{c\frac rn \log\log n}.\qquad\qquad\qquad
\eqno{(5.5)}
$$
Recall that $(q_0, q_1)=1$ and $\mathcal X_1 (\mod q_1)$ is primitive.

Expanding $F$ in additive characters $(\mod q_1)$ and applying the Gauss sum bound, it follows from the preceding that
$$
\Big|\sum_{v\equiv -u-w} F(v) \mathcal X_1(u+v+w)\Big|\ll (\log q_1) B^{C\frac rn\log\log n}\Big(1+\frac {2^{\Delta J}}{q_0q_1}\Big)\sqrt {q_1}.\eqno{(5.6)}
$$
Hence
$$
\begin{aligned}
(5.2)&< 2^{n-\Delta J-|A\backslash [J_1, J_2]|}. \ (5.6)\\
&< \frac {N 2^{-r}}{q_0} C^{\frac rn \Delta J} B^{C\frac rn \log\log n} \Big(q_0\sqrt{q_1} \, 2^{-\Delta J}+\frac 1{\sqrt{q_1}}\Big)\\
& <\frac {N2^{-r}} {q_0} \Big(\frac 1B +\frac 1{\sqrt{q_1}} \, B^{C\frac 4n\log\log n}\Big) \ \text { for appropriate $\Delta J\sim \log B$}\\
\end{aligned}
$$
$$
<\frac {N2^{-r}}{q_0} \frac 1{\sqrt{q_1}} \, B^{C\frac rn\log\log n}
\eqno{(5.7)}
$$
and it follows that
$$
\alpha(q_1, q_0)< 2^{-\frac n{12 r}}+\frac 1{\sqrt{q_1}} \, B^{C\frac rn \log\log n}.\eqno {(5.8)}
$$
Hence in (3.19), there is also the estimate
$$
\sum_{q_3<B} \frac {\alpha(q_1, q_3)}{q_3} < 2^{-\frac n{12r}} (\log B) +\frac{B^{C\frac rn\log\log n}}{\sqrt{\tilde q_1}} (\log B).\eqno{(5.10)}
$$
Combined with (4.31), it follows that for $\tilde q_1>1$
$$
\sum_{q_3<B} \ \frac {\alpha(q_1, q_3)}{q_3} < 2^{-\frac 13 10^{-3} \frac nr}\eqno{(5.11)}
$$
provided
$$
\log B< c \frac{n^2}{r^2}(\log\log n)^{-1}.
$$
as already implied by (4.24).

\section{Conclusion of the Argument}

In this section, we address the subdivision in `good' and `bad' characters and complete the argument.

Denote $N(\alpha, T;  \chi)$ the number of zero's of $L(s, \chi)$ such that
$$
\alpha\leq \sigma\leq 1, |t|\leq T \quad (s=\sigma+it).
$$
Then (see [Bom], Theorem 14)
$$
N(\alpha) =\sum_{q\leq Q} \quad \sideset{}{^*}\sum_{ \chi(\mod q)} \ N(\alpha, T; \chi)\ll (TQ)^{8(1-\alpha)}\eqno{(6.1)}
$$
where $\sum^*$ refers to summation over primitive characters.

Let $\chi$ be a non-principal character.
From Proposition 5.25 in [I-K], for $T\leq x$
$$
\psi (x,  \chi)= - \sum_{\substack {L(\rho,  \chi)=0\\ |\gamma|\leq T}} \ \frac {x^\rho -1}\rho +0\left(\frac xT(\log x q)^2\right)\eqno{(6.2)}
$$
where $\rho=\beta+ i\gamma$. We denote
$$
\eta =\eta(\chi) =\min (1-\beta)\eqno{(6.3)}
$$
with min taken over all zero's $\rho$ of $L(s, \chi)$ with $|\gamma|\leq T$.

The subdivision of characters in classes $\mathcal G$ and $\mathcal  B$ depends on whether $\eta \geq \eta^*$ (resp. $\eta<\eta^*)$ with $\eta^*$
to be specified.
Recall also that $q\leq  B$.

It follows from (6.2) that
$$
|\psi(x, \chi)|\leq x\left(\sum_\rho \frac 1{\rho x^{1-\beta}} + \frac {(\log x)^2}T\right)\eqno{(6.4)}
$$
and from the density bound and (6.3), assuming $qT< x^{\frac 1{20}}$
$$
\begin{aligned}
\sum \frac 1{x^{1-\beta}} &= -2 \int_{\frac 12}^{1-\eta} \frac 1{x^{1-\alpha}} dN(\alpha) = 2x^{-\frac 12} N\left(\frac 12\right) +2\log x
\int_\eta^{\frac 12} \left[\frac{(TQ)^8}x\right] ^\tau d\tau\\[6pt]
&<\frac{(qT)^4}{x^{1/2}} +x^{-\eta/2} < x^{-\eta/2}.
\end{aligned}
$$
Hence, for all $x$,
$$
|\psi(x, \chi)|< (QT)^{20}+\frac {x(\log x)^2}T+ x^{1-\frac \eta 2}.
\eqno{(6.5)}
$$
Fix $B<T=T(N)< N^{\frac 1{100}}$.
Assuming $\eta \geq \eta_*$, substitution of (6.5) in (2.7) gives $(r=|A|)$
$$
B (C\log n)^r \left\{(BT)^{20} +\frac {n^2N}T+ N^{1-\frac 12 \eta^* }\right\}.\eqno{(6.6)}
$$
If we let
$$
\eta_*=C \frac {\log T}{\log N}\eqno{(6.7)}
$$
(with appropriate $C>0$), the contribution of
$$
\mathcal G=\{\chi(\mod q); q<B, \eta(\chi)\geq \eta_*\}
$$
in (2.7) will be bounded by
$$
B(C\log n)^r \, n^2 \, \frac NT.\eqno{(6.8)}
$$
Next, consider the contribution of primitive $\chi_1(\mod q_1)$ in $\mathcal B$, given by (3.19).

Assuming $\tilde q_1>1$, (5.11) implies
$$
2^{-r}N(\log N)^5 \sum_{\chi_1\in \mathcal B} \, \frac {q_1}{\phi(q_1)} \sum_{\substack{q_3 < B\\ q_3 \, sf, \text  { odd }}}
\ \frac {\alpha(q_1, q_3)}{q_3} < 2^{-r} N n^6 \, 2^{-\frac 13 10^{-3}\frac nr} |\mathcal B|\eqno{(6.9)}
$$
where again by (6.1)
$$
|\mathcal  B| =|\{\chi_1 \text { primitive } (\mod q_1) \text { with } q_1<B \text
{ and } \eta(\chi_1)<\eta_*\}|<(TB)^{8\eta_*}.\eqno{(6.10)}
$$
Hence
$$
(6.9) < 2^{-r} N 2^{-\frac 14 10^{-3}\frac nr}\eqno{(6.11)}
$$
provided
$$
\log T< c\frac n{\sqrt r}\eqno{(6.12)}
$$
(for appropriate $c>0$).

Recall also assumption (4.24)
$$
\log B< c(\log\log n)^{-1} \frac {n^2}{r^2}.\eqno{(6.13)}
$$

Finally, consider the case $\tilde q_1 =1$, i.e. $q_1$ is a power of 2.

From the Gallagher-Iwaniec result (cf. [H-K], Lemma 5), we obtain the following improved zero-free region
$$
\eta(\chi_1) \gtrsim [(\log q_1 T)(\log\log q_1 T)]^{-3/4} \gg (\log T. \log\log T)^{-3/4}.\eqno{(6.14)}
$$
In order for $\chi$ to be induced by such $\chi_1$ to be in $\mathcal G$, we require
$$
\frac {\log T}{\log N} < c(\log T.\log\log T)^{-3/4}
$$
i.e
$$
\log T\ll n^{4/7}(\log n)^{-3/7}.\eqno{(6.15)}
$$

Collecting the bounds (1.11), (6.8), (3.15), (4.26), (3.19), (6.11), we obtain
$$
(1.1) - 2\mathbb E[f] N=
$$
$$
N2^{-r} O\left\{\frac {n^3(C\log n)^r}{\sqrt B} +\frac {n^2 B (C\log n)^r}T+B \, 2^r \exp \left(-c\frac {n^{3/5}}{(\log n)^{1/5}}\right)
+2^{-10^{-4}\frac nr}\right\}.\eqno{(6.16)}
$$
Fix $0<\gamma<1$. If we let
$$
B=n^6 (C\log n)^{2r}\gamma^{-2} \text { and }T=n^{8} (C\log n)^{3r} \gamma^{-3}
$$
and in order to fulfill (6.12), (6.13), (6.15), we assume
$$
r < \left(\frac n{\log n}\right)^{4/7}\eqno{(6.17)}
$$
then (6.16) implies
$$
(1.1)- 2\mathbb E[f]N<\gamma \mathbb E[f].
$$
This proves the following

\begin{theorem}
{\it Let $N=2^n$, $n$ large enough.
Under assumption (6.17), choosing a set $A\subset \{1, \ldots, n\}$, $|A|=r$ of binary digits and assigned values $\alpha_j=0, 1$ for $j\in A$, we have
$$
|\{p<N; p \text { has $j$-digit equal to $\alpha_j$ for $j\in A$} \}|=
$$
$$
\big(1+o(1)\big) \frac N{\log N} 2^{-r}.
$$}
\end{theorem}

The argument permits to give a better error term that we do not specify here.

\noindent
{\bf Acknowledgement.}  The author is grateful to the referee for his comments and suggestions to improve the presentation.
\null\vskip.2 true in

\end{document}